\documentclass[aps,pra,amssymb,amsmath,eqsecnum,nofootinbib]{revtex4}
\usepackage{graphicx}

\newtheorem{definition}{Definition}[section]

\newtheorem{proposition}{Proposition}[section]

\newtheorem{corollary}{Corollary}[section]

\newtheorem{remark}{Remark}[section]
\newtheorem{axiom}{AXIOM}[section]

\newenvironment{hypothesis}{HP: \begin{center}} {\end{center}}
\newenvironment{thesis}{TH: \begin{center}} {\end{center}}
\newenvironment{proof}{\begin{center}PROOF: \end{center}} {$ \blacksquare $}

\begin{document}
\title{A new kind of numbers, the Non-Dedekindian Numbers, and the extension to them of the notion of algorithmic randomness}
\author{Gavriel Segre}
\date{8 December 2006}
\begin{abstract}
  A new number system, the set of the non-Dedekindian numbers, is introduced and characterized axiomatically.

It is then proved that any hypercontinous hyperreal number system
is strictly included in the set of the Non-Dedekindian Numbers.

The notion of algorithmic-randomness is then extended to
non-Dedekindian numbers.

As a particular case, the notion of algorithmic randomness for the
particular hyperreal number system of Non-Standard Analysis is
explicitly analyzed.
\end{abstract}
\maketitle \footnote{The reported date is the one of the first
public (i.e. appeared on my homepage http://www.gavrielsegre.com)
version of this paper that has, anyway, received later
improvements.}
\newpage
\tableofcontents
\newpage
\section{Introduction}

The non-euclidean revolution \cite{Trudeau-87} consisted in the
discovery of many kind of non-euclidean geometries corresponding
to different choices in the imposed cardinality of the set of
straight lines parallel to a given straight line and passing for a
point not belonging to it (euclidean geometry consisting in the
assumption, stated by Euclides' Fifth Axiom, that such a cardinal
number is equal to one).

In a completely different framework, the axiomatic definition of
the set $ \mathbb{R} $ of the real numbers \cite{Schechter-98},
Dedekind's Continuity Axiom resembles, with this respect,
Euclides' Fifth Axiom since again it imposes that the cardinality
of a suitable set (the intersection of a sequence of nested
halving intervals) is equal to one.

Such a similarity naturally induces to investigate which kind of
different numbers' systems we obtain by replacing Dedekind's Axiom
with  different choices in the imposed cardinality of the
intersection of nested halving intervals.

\smallskip

What one obtains in this way is a new number system, the set $
\mathbb{ND} $ of the Non-Dedekindian numbers, that here we
explicitly define axiomatically.

\smallskip

We then show that any hypercontinuous hyperreal number system is
strictly included in $ \mathbb{ND} $.

\smallskip

Finally we extend the notion of algorithmic-randomness to
non-Dedekindian numbers.

As a particular case, the notion of algorithmic randomness for the
particular hyperreal number system of Non-Standard Analysis is
explicitly analyzed.

\newpage
\section{Different ways to define the real numbers}

As it is well known there exist many different equivalent ways of
defining the set $ \mathbb{R} $ of the real numbers.

Many of them has the same structure: they define $ \mathbb{R} $ as
a chain-ordered field satisfying some supplementary condition of
completeness (belonging to a family of equivalent completeness'
conditions).

The simpler approach was given by Cantor: introduced on $
\mathbb{Q} $ the metric $ d ( r_{1} , r_{2} ) := | r_{1} - r_{2} |
$ one defines $ \mathbb{R} $ as the metric completion of the
metric space $ ( \mathbb{Q} , d ) $.

In this way a real number is then defined as an equivalence class
of Cauchy sequences over $  \mathbb{Q} $ with respect to the
following equivalence relation:
\begin{equation}
    \{ r_{n}\}_{n \in \mathbb{N}} \; \sim \;   \{ s_{n}\}_{n \in
    \mathbb{N}} \; := \; \lim_{n \rightarrow + \infty} d ( r_{n} ,
    s_{n} ) \, = \, 0
\end{equation}

In this paper we will concentrate, anyway, our attention on
Dedekind's way of formalizing the supplementary completeness
condition since:
\begin{enumerate}
    \item it shows in such an intuitive way the evidence that  the
    supplementary added completeness condition is a condition assuring the
    "continuity" of $ \mathbb{R} $ that it is usually called the
    Continuity Axiom
    \item it has a natural link with Algorithmic Information
    Theory
    \item it has a natural generalization that will allow us to
    introduce a new number system: the set $ \mathbb{ND}$ of the non-Dedekindian
    numbers
\end{enumerate}

\newpage
\section{Dedekind Continuity Axiom}

Let us denote by $ \Sigma := \{ 0 , 1 \} $ the binary alphabet, by
$ \Sigma^{\star} := \cup_{n \in \mathbb{N}_{+}} \Sigma^{n} $ the
set of all the binary strings and by $ \Sigma^{\infty} := \{
\bar{x} : \mathbb{N}_{+} \mapsto \Sigma \} $ the set of all the
binary sequences.  Given $ \bar{x} \in \Sigma^{\infty} $ and $
a_{0} , b_{0} \in \mathbb{R} : a_{0} < b_{0} $ let us introduce
the following definition by induction:
\begin{definition}
\end{definition}
\begin{itemize}
    \item
\begin{equation}
    a_{1} ( \bar{x} ) \; := \; \left\{%
\begin{array}{ll}
    a_{0} & \hbox{if $  x_{1} = 0$} \\
    \frac{a_{0}+b_{0}}{2} & \hbox{ if $ x_{1} = 1$} \\
\end{array}%
\right.
\end{equation}
\begin{equation}
    b_{1} ( \bar{x} ) \; := \; \left\{%
\begin{array}{ll}
     \frac{a_{0}+b_{0}}{2}  & \hbox{ if $ x_{1} = 0$} \\
    b_{0}  & \hbox{ if $ x_{1} = 1$} \\
\end{array}%
\right.
\end{equation}
    \item
\begin{equation}
    a_{n} ( \bar{x} ) \; := \; \left\{%
\begin{array}{ll}
    a_{n-1} & \hbox{if $ x_{n} = 0$} \\
    \frac{a_{n-1}+b_{n-1}}{2} & \hbox{if $ x_{n} = 1$} \\
\end{array}%
\right.
\end{equation}
\begin{equation}
    b_{n} ( \bar{x} ) \; := \; \left\{%
\begin{array}{ll}
     \frac{a_{n-1}+b_{n-1}}{2}  & \hbox{if $  x_{n} = 0$} \\
    b_{n-1} & \hbox{if $ x_{n} = 1$} \\
\end{array}%
\right.
\end{equation}
\end{itemize}

We can at last introduce the following:
\begin{definition} \label{def:real numbers}
\end{definition}
\emph{set $ \mathbb{R}$ of the real numbers:}

a chain-ordered field satisfying the following axiom
\ref{ax:Dedekind Continuity Axiom}

\begin{axiom} \label{ax:Dedekind Continuity Axiom}
\end{axiom}
\emph{Dedekind Continuity Axiom:}
\begin{equation}
    \exists  \,!  \, N_{Dedekind} ( a_{0} , b_{0} ; \bar{x} ) \in \mathbb{R} \; : \; N_{Dedekind} ( a_{0} , b_{0} ; \bar{x} ) \in  \cap_{n \in \mathbb{N}} [ a_{n} ( \bar{x} ) , b_{n}
    ( \bar{x} ) ] \ \; \; \forall a_{0} , b_{0} \in \mathbb{R} : a_{0} < b_{0} \, , \,  \forall \bar{x} \in \Sigma^{\infty}
\end{equation}

It may be proved that:

\begin{proposition} \label{prop:equivalence between Dedekind Continuity Axiom and Dedekind completeness}
\end{proposition}

\begin{hypothesis}
\end{hypothesis}

\begin{center}
  $ ( F , \leq ) $ chain-ordered field
\end{center}

\begin{thesis}
\end{thesis}
\begin{center}
  Axiom \ref{ax:Dedekind Continuity Axiom} is equivalent to
  Dedekind completeness
\end{center}

\smallskip

\begin{remark}
\end{remark}
Owing to Proposition \ref{prop:equivalence between Dedekind
Continuity Axiom and Dedekind completeness} the definition
\ref{def:real numbers} of $ \mathbb{R} $ is equivalent to the more
usual definition of $ \mathbb{R} $ as a Dedekind complete
chain-ordered field.

Axiom \ref{ax:Dedekind Continuity Axiom} has, anyway, a
constructive nature that lacks to the condition of Dedekind
Completeness:

given $ n \in \mathbb{N} $, $ a_{n} ( \bar{x} ) $ and $ b_{n} (
\bar{x} ) $ may be computed through the following Mathematica
\cite{Wolfram-96} expressions:
\begin{verbatim}

leftextreme[n_, a_, b_, binarystring_] :=
    If[n == 1, If[Part[binarystring, 1] == 0, a, (a + b)/2] ,
      If[Part[binarystring, n] == 0, leftextreme[n - 1, a, b, binarystring],
        (leftextreme[n - 1, a, b, binarystring] +
            rightextreme[n - 1, a, b, binarystring])/2]]


rightextreme[n_, a_, b_, binarystring_] :=
    If[n == 1, If[Part[binarystring, 1] == 0, (a + b)/2, b],
      If[Part[binarystring, n] == 0,
        (leftextreme[n - 1, a, b, binarystring] +
            rightextreme[n - 1, a, b, binarystring])/2,
        rightextreme[n - 1, a, b, binarystring]]]

\end{verbatim}

\bigskip

Let us now introduce the following:
\begin{definition}
\end{definition}
\emph{Dedekind operator with respect to $ [ a_{0} , b_{0}]$:}

$ \mathcal{D}_{a_{0},b_{0}} : \Sigma^{\infty} \mapsto [ a_{0} ,
b_{0}]  $:
\begin{equation}
     \mathcal{D}_{a_{0},b_{0}} ( \bar{x} )  \; := \;  N_{Dedekind} ( a_{0} , b_{0} ; \bar{x} )
\end{equation}

The axiom \ref{ax:Dedekind Continuity Axiom} implies that:
\begin{corollary}
\end{corollary}
\begin{equation}
    \mathcal{D}_{a_{0},b_{0}}  \text{ is bijective } \; \; \forall
    a_{0} , b_{0} \in \mathbb{R} \; : \; a_{0} < b_{0}
\end{equation}
\begin{proof}
\begin{enumerate}
    \item Let us prove that $  \mathcal{D}_{a_{0},b_{0}} $ is
    injective.

Given $ \bar{x}, \bar{y} \in \Sigma^{\infty} :  \bar{x} \neq
\bar{y} $ this means that:
\begin{equation}
    \exists n \in \mathbb{N} : x_{n} \neq y_{n}
\end{equation}
But then:
\begin{equation}
    [ a_{n} ( \bar{x} ) , b_{n}  ( \bar{x} ) ] \neq   [ a_{n} ( \bar{y} ) , b_{n}  ( \bar{y} ) ]
\end{equation}
and hence:
\begin{equation}
     \mathcal{D}_{a_{0},b_{0}}  ( \bar{x} ) \; \neq \; \mathcal{D}_{a_{0},b_{0}}  ( \bar{y} )
\end{equation}
    \item  Let us prove that $  \mathcal{D}_{a_{0},b_{0}} $ is
    surjective.

Given $ c \in [ a_{0} , b_{0}] $ let us choose at each step the
value of $ x_{n} \in \Sigma $ such that $ c \in [ a_{n} , b_{n} ]
$.

Then by construction:
\begin{equation}
     \mathcal{D}_{a_{0},b_{0}}  ( \bar{x} ) \; = \; c
\end{equation}
\end{enumerate}
\end{proof}

\newpage
\section{Algorithmic information theoretic analysis of the Dedekind operator}
Let us observe, first of all, that denoted with $ v_{2}$  the
2-ary map introduced in the definition \ref{def:n-ary value} of
section \ref{sec:algorithmic randomness}:

\begin{proposition}
\end{proposition}
\begin{equation}
     \mathcal{D}_{0,1} ( \bar{x} )  \; = \;
   v_{2} ( \bar{x} ) \; \; \forall \bar{x}
    \in \Sigma^{\infty}
\end{equation}
\begin{proof}
By construction:
\begin{equation}
    \lim_{n \rightarrow + \infty} a_{n} ( \bar{x} ) \; = \;
    \sum_{n=1}^{\infty}\frac{x_{n}}{2^{n}} \; \; \forall \bar{x}
    \in \Sigma^{\infty}
\end{equation}
But axiom \ref{ax:Dedekind Continuity Axiom} implies that:
\begin{equation}
    \lim_{n \rightarrow + \infty} a_{n} ( \bar{x} ) \; = \; \lim_{n \rightarrow + \infty} b_{n} ( \bar{x}
    ) \; = \;   N_{Dedekind} ( 0 , 1 ; \bar{x} ) \; \; \forall \bar{x}
    \in \Sigma^{\infty}
\end{equation}
\end{proof}

i.e. $ N_{Dedekind} ( 0 , 1 ; \bar{x} ) $ is the number having $
(0.\bar{x})_{2}  $ as base-two representation.

Hence \cite{Calude-02}:
\begin{corollary}
\end{corollary}
\begin{enumerate}
    \item
\begin{equation}
     \mathcal{D}_{0,1} ( \bar{x} ) \in \mathbb{Q} \; \Leftrightarrow \; \exists
     \vec{y} \in \Sigma^{\star} \cup \{ \lambda \} , \vec{z} \in
     \Sigma^{\star} \; : \; \bar{x} = \vec{y} \cdot \vec{z}^{ \, \infty}
\end{equation}
where $ \cdot $ denotes concatenation, $ \lambda$ is  the empty
string and where $ \vec{z}^{ \, \infty} $ denotes the infinite
repetition of the string $\vec{z}$.
    \item
\begin{equation}
  \mathcal{D}_{0,1} (RANDOM ( \Sigma^{\infty} ))  \; = \;  RANDOM([ 0 , 1])
\end{equation}
where $  RANDOM( 0 , 1) $ is the set of the random reals belonging
to the interval $ [0 ,1] $ while $ RANDOM ( \Sigma^{\infty} )$ is
the set of the random binary sequences.
\end{enumerate}

Let us now consider two arbitrary $ a_{0} , b_{0} \in \mathbb{R} :
a_{0} < b_{0} $ and let us introduce the following set:

\begin{definition} \label{def:algorithmically random numbers of  an interval}
\end{definition}
\emph{real random numbers with respect to $ [a_{0} , b_{0}] $:}
\begin{equation}
    RANDOM ( [a_{0} , b_{0}] ) \; := \;  \mathcal{D}_{a_{0},b_{0}} (RANDOM ( \Sigma^{\infty} ))
\end{equation}

Let us introduce the following map:
\begin{definition}
\end{definition}
$\mathcal{T}_{a_{0} , b_{0}} : [ a_{0} , b_{0} ]  \mapsto [ 0 ,
1]$:
\begin{equation}
    \mathcal{T}_{a_{0} , b_{0}}  (x) \; := \; \frac{x - a_{0}}{b_{0} - a_{0}}
\end{equation}
and its inverse:
\begin{definition}
\end{definition}
$\mathcal{T}_{a_{0} , b_{0}}^{-1} : [ 0 , 1] \mapsto [ a_{0} ,
b_{0} ]$:
\begin{equation}
  \mathcal{T}_{a_{0} , b_{0}}^{-1} (x) \; := \; ( b_{0} - a_{0})x
  + a_{0}
\end{equation}

Then by construction:
\begin{proposition}
\end{proposition}
\begin{equation}
     \mathcal{T}_{a_{0} , b_{0}} (  \mathcal{D}_{a_{0} , b_{0}}   ( \bar{x}
     )) \; = \; v_{2} ( \bar{x} ) \; \;
     \forall \bar{x} \in \Sigma^{\infty}
\end{equation}

\smallskip

\begin{corollary}
\end{corollary}
\begin{enumerate}
    \item
\begin{equation}
    \mathcal{T}_{a_{0} , b_{0}} ( \mathcal{D}_{a_{0},b_{0}} ( \bar{x} )) \in \mathbb{Q} \; \Leftrightarrow \; \exists
     \vec{y} \in \Sigma^{\star} \cup \{ \lambda \} , \vec{z} \in
     \Sigma^{\star} \; : \; \bar{x} = \vec{y} \cdot \vec{z}^{ \, \infty}
\end{equation}
    \item
\begin{equation}
  \mathcal{T}_{a_{0} , b_{0}} (RANDOM [ a_{0} , b_{0} ] )  \; = \;  RANDOM( 0 , 1)
\end{equation}
\end{enumerate}

\smallskip

\begin{remark}
\end{remark}
It is extremely important to remark, at this point , that while $
RANDOM( 0 , 1) $ , according to definition \ref{def:random numbers
between zero and one}, is an intrinsic notion  characterizing the
set of the random reals belonging to the interval $ [ 0,1 ] $, $
RANDOM ([ a_{0} , b_{0} ]) $ is not an intrinsic notion that
characterizes the set of the random numbers belonging to $ [a_{0}
, b_{0} ] $ but a relative notion that characterizes the random
reals with respect to $  [a_{0} , b_{0} ] $, i.e. the set of the
random reals in the interval $[ 0,1]$ when such an interval is
seen dilatated and translated by $ \mathcal{T}_{a_{0},b_{0}}^{- 1}
$.

That this is the case may be appreciated considering the
following:
\begin{proposition}
\end{proposition}
\begin{hypothesis}
\end{hypothesis}
\begin{equation}
    a_{1}, a_{2}, b_{1}, b_{2} \in \mathbb{R} \; : \; [ a_{2} ,
    b_{2}] \, \subset \, [ a_{1} ,
    b_{1}]
\end{equation}

\begin{thesis}
\end{thesis}
\begin{equation}
    RANDOM( [ a_{2} , b_{2} ] ) \; \neq \;  RANDOM( [ a_{1} , b_{1} ]
    ) \cap [ a_{2} , b_{2} ]
\end{equation}
\begin{proof}
Since:
\begin{equation}
    RANDOM(0,1) \; = \; \mathcal{T}_{a_{1},b_{1}}^{- 1} \{ RANDOM ( [ a_{1} , b_{1} ] )
    \} \; = \;  \mathcal{T}_{a_{2},b_{2}}^{- 1} \{ RANDOM ( [ a_{2} , b_{2} ] )
    \}
\end{equation}
if follows that:
\begin{equation}
    RANDOM ( [ a_{2} , b_{2} ]) \; = \; \frac{( b_{1} - a_{1}) RANDOM ( [ a_{1} , b_{1} ] ) + a_{1} -
    a_{2}}{b_{2} -  a_{2}} \; \neq \; RANDOM( [ a_{1} , b_{1} ]
    ) \cap [ a_{2} , b_{2} ]
\end{equation}
\end{proof}
\newpage
\section{A new number system: the Non-Dedekindian numbers}

Given $ n \in \mathbb{N}_{+} \cup \{ \aleph_{n} , n \in \mathbb{N}
\} $  \footnote{We have assumed that  $ n < \aleph_{\omega}$. This
has been done since assumed the following:
\begin{axiom} \label{ax:Generalized Continuum Hypothesis}
\end{axiom}
\emph{Generalized Continuum Hypothesis:}
\begin{equation}
    2^{\aleph_{n}} \; = \; \aleph_{n+1} \; \; \forall n \in \mathbb{N}
\end{equation}
it follows that:
\begin{equation}
   \aleph_{n} \; = \; | \mathcal{P}^{n} ( \mathbb{N} ) | \; \;
   \forall n \in \mathbb{N}
\end{equation}
(where $ | S| $ denotes the cardinality of a set S and where $
\mathcal{P}^{n} $ denotes the $ n^{th} $ iterate of the power-set
operator) while cardinals $ \geq \aleph_{\omega} $ cannot be
obtained in this way.} let $ \mathbb{G}_{n} $ be a chain-ordered
field \cite{Schechter-98}.

 Given $ \bar{x} \in \Sigma^{\infty} $
and $ a_{0} , b_{0} \in \mathbb{G}_{n} : a_{0} < b_{0} $ let us
introduce  the following definition by induction:
\begin{definition}
\end{definition}
\begin{itemize}
    \item
\begin{equation}
    a_{1} ( \bar{x} ) \; := \; \left\{%
\begin{array}{ll}
    a_{0} & \hbox{if $  x_{1} = 0$} \\
    \frac{a_{0}+b_{0}}{2} & \hbox{ if $ x_{1} = 1$} \\
\end{array}%
\right.
\end{equation}
\begin{equation}
    b_{1} ( \bar{x} ) \; := \; \left\{%
\begin{array}{ll}
     \frac{a_{0}+b_{0}}{2}  & \hbox{ if $ x_{1} = 0$} \\
    b_{0}  & \hbox{ if $ x_{1} = 1$} \\
\end{array}%
\right.
\end{equation}
    \item
\begin{equation}
    a_{n} ( \bar{x} ) \; := \; \left\{%
\begin{array}{ll}
    a_{n-1} & \hbox{if $ x_{n} = 0$} \\
    \frac{a_{n-1}+b_{n-1}}{2} & \hbox{if $ x_{n} = 1$} \\
\end{array}%
\right.
\end{equation}
\begin{equation}
    b_{n} ( \bar{x} ) \; := \; \left\{%
\begin{array}{ll}
     \frac{a_{n-1}+b_{n-1}}{2}  & \hbox{if $  x_{n} = 0$} \\
    b_{n-1} & \hbox{if $ x_{n} = 1$} \\
\end{array}%
\right.
\end{equation}
\end{itemize}

\begin{definition}
\end{definition}
\emph{Dedekind set with respect to $ [ a_{0} , b_{0}] $ and
$\bar{x}$:}
\begin{equation}
    S_{Dedekind} ( a_{0} , b_{0} ; \bar{x} ) \; := \; \{ d \in
    \mathbb{G}_{n} \, : d \, \in  \cap_{n \in \mathbb{N}} [ a_{n} ( \bar{x} ) , b_{n}
    ( \bar{x} ) ] \}
\end{equation}

Let us now introduce the following:
\begin{definition}
\end{definition}
\emph{set of the generalized numbers of order n:}
\begin{center}
  the chain-ordered field $ \mathbb{G}_{n} $ satisfying
  the following axiom \ref{ax:generalized Dedekind Axiom of order n}
\end{center}

\newpage

\begin{axiom} \label{ax:generalized Dedekind Axiom of order n}
\end{axiom}
\emph{Generalized Dedekind Axiom of order n:}
\begin{equation}
    |  S_{Dedekind} ( a_{0} , b_{0} ; \bar{x} ) | \; = \; n \; \;
    \forall a_{0} , b_{0} \in \mathbb{G}_{n} : a_{0} < b_{0} \, ,
    \, \forall \bar{x} \in \Sigma^{\infty}
\end{equation}
where $ | S | $ denotes the cardinality of a set S.

 Clearly:
\begin{proposition} \label{prop:generalized number of first order are the reals}
\end{proposition}
\begin{equation}
    \mathbb{G}_{1} \; = \; \mathbb{R}
\end{equation}
\begin{proof}
 For $ n=1 $ the axiom \ref{ax:generalized Dedekind Axiom of order
 n} reduces to Dedekind's Continuity Axiom.
\end{proof}

\smallskip

Given $ n \in  \mathbb{N}_{+} \cup \{ \aleph_{n} , n \in
\mathbb{N} \}
 \, : \, n > 1 $ we will call $ \mathbb{G}_{n} $ the \emph{set of the
non-Dedekindian numbers of order n} and we will call $ \mathbb{ND}
 \, := \, \cup_{ \{ n \in \mathbb{N}_{+} : n > 1 \} \cup \{ \aleph_{n} , n \in \mathbb{N}
\}}  \mathbb{G}_{n} $ the set of the \emph{Non-Dedekindian
numbers}.

\smallskip

\begin{remark}
\end{remark}
Clearly the furnished formal axiomatic definition of the
Non-Dedekindian numbers is not the whole story. One has:
\begin{enumerate}
    \item to prove that the involved formal system is consistent
    \item to prove that $ \mathbb{ND} \; \neq \; \emptyset $
\end{enumerate}

\smallskip

\begin{remark}
\end{remark}
As to the proof of the consistence of the given axiomatic
definition of generalized numbers of any order, let us remark that
we know that for $ n=1$ this is true since we know that Dedekind's
axiomatization of real numbers is consistent.

Hence, to obtain a proof by induction, it would be sufficient to
prove that the consistence of the axioms for $ \mathbb{G}_{n} $
implies the consistence for the axioms for $ \mathbb{G}_{n+1} $.

We leave this task open for future investigation.

\smallskip

The proof that $ \mathbb{ND} \; \neq \; \emptyset $ will be given
in the next section where we will prove that any hypercontinuous
hyperreal number system is strictly included in $ \mathbb{ND}$.

\newpage
\section{Hypercontinuous Hyperreal Numbers as particular Non Dedekindian Numbers}
Let us start from the following:
\begin{definition} \label{def:hyperreal number system}
\end{definition}
\emph{hyperreal number system:}
\begin{center}
  a chain-ordered non-Archimedean field containing $ \mathbb{R}$ as a subfield.
\end{center}

Given an  hyperreal number system $ \mathbb{H} $:
\begin{definition}
\end{definition}
\emph{non-standard part of $ \mathbb{H} $:}
\begin{equation}
  \mathbb{H}_{non standard} \; := \; \mathbb{H} \, - \, \mathbb{R}
\end{equation}
\begin{definition}
\end{definition}
\emph{infinitesimals of $ \mathbb{H}$:}
\begin{equation}
  \mathbb{H}_{infinitesimals} \; := \; \{ x \in \mathbb{H} \, : \, x \in ( -r , r ) \; \forall r \in \mathbb{R}_{+}   \}
\end{equation}
\begin{definition}
\end{definition}
\emph{unboundeds of $ \mathbb{H}$:}
\begin{equation}
  \mathbb{H}_{unboundeds} \; := \; \{ x \in \mathbb{H} \, : \,
  \nexists \, r \in \mathbb{R}_{+} : x \in ( -r , r ) \}
\end{equation}

\smallskip

It may be proved that \cite{Schechter-98}:

\begin{proposition}
\end{proposition}
\begin{enumerate}
    \item
\begin{equation}
  \mathbb{H}_{infinitesimals} \cap \mathbb{R} \; = \;  \{ 0 \}
\end{equation}
    \item
\begin{equation}
 x \in \mathbb{H}_{infinitesimals} - \{ 0 \} \; \Leftrightarrow \;
 \frac{1}{x} \in  \mathbb{H}_{unboundeds}
\end{equation}
 \item
 \begin{equation}
    |   \mathbb{H}_{infinitesimals} | \; \geq \; \aleph_{0}
\end{equation}

\end{enumerate}

so that obviously:
\begin{corollary}
\end{corollary}
\begin{equation}
    | \mathbb{H}_{unboundeds} | \; = \;  |  \mathbb{H}_{infinitesimals}
    | \; \geq \; \aleph_{0}
\end{equation}

Let us introduce the following:

\begin{definition} \label{def:hypercontinuity}
\end{definition}
\emph{$ \mathbb{H}$ is hypercontinuous:}
\begin{enumerate}
    \item
\begin{equation}
    | \mathbb{H}_{infinitesimals} | \; \geq \; \aleph_{1}
\end{equation}
    \item
\begin{equation}
    | [ \epsilon_{1} , \epsilon_{2} ] | \; = \; |
    \mathbb{H}_{infinitesimals} | \; \; \forall \epsilon_{1},
    \epsilon_{2} \in  \mathbb{H}_{infinitesimals} \; : \;  \epsilon_{1} \, < \, \epsilon_{2}
\end{equation}
\end{enumerate}

\smallskip

\begin{remark}
\end{remark}
We would like to caution the reader that the terminology of
definition \ref{def:hypercontinuity} is new.

\smallskip

\smallskip

Given  $ x_{1} , x_{2} \in \mathbb{H} $ let us introduce the
following:
\begin{definition}
\end{definition}
\emph{$x_{1}$ is infinitesimally closed to $x_{2}$: }
\begin{equation}
    x_{1} \sim_{\text{infinitesimally closed}} x_{2} \; := \;
    x_{1} - x_{2} \in  \mathbb{H}_{infinitesimals}
\end{equation}
and let us recall that \cite{Schechter-98}:
\begin{proposition}
\end{proposition}
\begin{enumerate}
    \item $ \sim_{\text{infinitesimally closed}} $ is an
    equivalence relation over $ \mathbb{H} $.
    \item
\begin{equation}
    \forall x \in \mathbb{R} \cup \mathbb{H}_{infinitesimals} \; \exists \, ! \, std (x) \in
    \mathbb{R} \, : \; x \sim_{\text{infinitesimally closed}}
    std(x)
\end{equation}
\end{enumerate}
$std(x)$ is called the \emph{standard part} of $ x \in
\mathbb{H}$.

\smallskip

Let us now repeat the construction of the previous section.

Given $ \bar{x} \in \Sigma^{\infty} $ and $ a_{0} , b_{0} \in
\mathbb{H} : a_{0} < b_{0} $ let us introduce  the following
definition by induction:
\begin{definition}
\end{definition}
\begin{itemize}
    \item
\begin{equation}
    a_{1} ( \bar{x} ) \; := \; \left\{%
\begin{array}{ll}
    a_{0} & \hbox{if $  x_{1} = 0$} \\
    \frac{a_{0}+b_{0}}{2} & \hbox{ if $ x_{1} = 1$} \\
\end{array}%
\right.
\end{equation}
\begin{equation}
    b_{1} ( \bar{x} ) \; := \; \left\{%
\begin{array}{ll}
     \frac{a_{0}+b_{0}}{2}  & \hbox{ if $ x_{1} = 0$} \\
    b_{0}  & \hbox{ if $ x_{1} = 1$} \\
\end{array}%
\right.
\end{equation}
    \item
\begin{equation}
    a_{n} ( \bar{x} ) \; := \; \left\{%
\begin{array}{ll}
    a_{n-1} & \hbox{if $ x_{n} = 0$} \\
    \frac{a_{n-1}+b_{n-1}}{2} & \hbox{if $ x_{n} = 1$} \\
\end{array}%
\right.
\end{equation}
\begin{equation}
    b_{n} ( \bar{x} ) \; := \; \left\{%
\begin{array}{ll}
     \frac{a_{n-1}+b_{n-1}}{2}  & \hbox{if $  x_{n} = 0$} \\
    b_{n-1} & \hbox{if $ x_{n} = 1$} \\
\end{array}%
\right.
\end{equation}
\end{itemize}

\begin{definition}
\end{definition}
\emph{Dedekind set with respect to $ [ a_{0} , b_{0}] $ and
$\bar{x}$:}
\begin{equation}
    S_{Dedekind} ( a_{0} , b_{0} ; \bar{x} ) \; := \; \{ d \in
    \mathbb{H} \, : d \, \in  \cap_{n \in \mathbb{N}} [ a_{n} ( \bar{x} ) , b_{n}
    ( \bar{x} ) ] \}
\end{equation}

\smallskip

Then:

\newpage
\begin{proposition} \label{prop:first case}
\end{proposition}

\begin{hypothesis}
\end{hypothesis}
\begin{equation}
  a_{0} , b_{0} \in \mathbb{R} \cup  \mathbb{H}_{infinitesimals}
  \; \wedge \; std ( a_{0} ) \neq std ( b_{0} )
\end{equation}
\begin{thesis}
\end{thesis}
\begin{enumerate}
    \item
\begin{equation}
| \{ std (d) , d \in  S_{Dedekind} ( a_{0} , b_{0} ; \bar{x} ) \}
| \; = \; 1 \; \; \forall \bar{x} \in \Sigma^{\infty}
\end{equation}
    \item
\begin{equation}
   | S_{Dedekind} ( a_{0} , b_{0} ; \bar{x} ) | \; = \; |
   \mathbb{H}_{infinitesimals} | \; \;  \forall \bar{x} \in \Sigma^{\infty}
\end{equation}
\end{enumerate}
\begin{proof}
\begin{enumerate}
    \item Let us observe first of all that:
\begin{equation}
    std( a_{n} ( \bar{x} ) )\; \neq \; std( b_{n}  ( \bar{x} )) \; \; \forall n \in
    \mathbb{N} , \forall \bar{x} \in \Sigma^{\infty}
\end{equation}
The thesis follows applying the axiom \ref{ax:Dedekind Continuity
    Axiom} to the set $ [ std(a_{0}) , std(b_{0}) ] \cap
    \mathbb{R} $
    \item the thesis follows observing that:
\begin{equation}
    d + \epsilon \; \in [ a_{n} ( \bar{x} ) , b_{n} ( \bar{x} ) ] \; \; \forall n \in
    \mathbb{N}, \forall \epsilon
    \in \mathbb{H}_{infinitesimals} ,  \forall \bar{x} \in \Sigma^{\infty}
\end{equation}
\end{enumerate}
\end{proof}

\smallskip

\begin{proposition} \label{prop:second case}
\end{proposition}
\begin{hypothesis}
\end{hypothesis}
\begin{enumerate}
    \item $ \mathbb{H} $ is hypercontinuous
    \item
\begin{equation}
  a_{0} , b_{0} \in \mathbb{R} \cup  \mathbb{H}_{infinitesimals}
  \; \wedge \; std ( a_{0} ) = std ( b_{0} )
\end{equation}
\end{enumerate}

\begin{thesis}
\end{thesis}
 \begin{equation}
     | S_{Dedekind} ( a_{0} , b_{0} ; \bar{x} ) | \; = \; |
   \mathbb{H}_{infinitesimals} | \; \;  \forall \bar{x} \in \Sigma^{\infty}
\end{equation}
\begin{proof}
The thesis follows applying the definition
\ref{def:hypercontinuity}.
\end{proof}

\smallskip

Let us now include unbounded hyperreals in the game.

Given $ x_{1} , x_{2} \in \mathbb{H} $:
\newpage
\begin{definition}
\end{definition}
\emph{$x_{1}$ is finitely distant from  $x_{2}$: }
\begin{equation}
    x_{1} \, \sim_{\text{finitely distant from}} \, x_{2} \; :=
    \: x_{1}-x_{2} \in \mathbb{R} \cup \mathbb{H}_{infinitesimals}
\end{equation}

Then \cite{Schechter-98}:
\begin{proposition}
\end{proposition}
\begin{center}
 $ \sim_{\text{finitely distant from}} $ is an equivalence
 relation over $ \mathbb{H} $
\end{center}

Let us now consider the various cases as to the cardinality of the
Dedekind set:
\begin{proposition}
\end{proposition}
\begin{hypothesis}
\end{hypothesis}
\begin{enumerate}
    \item $ \mathbb{H} $ is hypercontinuous
    \item
\begin{equation}
     a_{0} \sim_{\text{finitely distant from}} b_{0}
\end{equation}
\end{enumerate}

\begin{thesis}
\end{thesis}
\begin{equation}
    | S_{Dedekind} ( a_{0} , b_{0} ; \bar{x} ) | \; = \; |
    \mathbb{H}_{infinitesimals} | \; \; \forall \bar{x} \in \Sigma^{\infty}
\end{equation}
\begin{proof}
  Clearly:
\begin{equation}
     a_{n} ( \bar{x} )  \sim_{\text{finitely distant from}}  b_{n} ( \bar{x} ) \; \; \forall n \in
     \mathbb{N} , \forall  \bar{x} \in \Sigma^{\infty}
\end{equation}
the thesis follow applying Proposition \ref{prop:first case} and
Proposition \ref{prop:second case}.
\end{proof}

\begin{proposition} \label{prop:third case}
\end{proposition}
\begin{hypothesis}
\end{hypothesis}
\begin{enumerate}
    \item $ \mathbb{H} $ is hypercontinuous
    \item
\begin{equation}
 a_{0} \nsim_{\text{finitely distant from}} b_{0}
\end{equation}
\end{enumerate}
\begin{thesis}
\end{thesis}
\begin{equation}
    | S_{Dedekind} ( a_{0} , b_{0} ; \bar{x} ) | \; = \; |
    \mathbb{H}_{infinitesimals} | \; \; \forall \bar{x} \in \Sigma^{\infty}
\end{equation}
\newpage
\begin{proof}
Clearly:
\begin{equation}
     a_{n} ( \bar{x} )  \nsim_{\text{finitely distant from}}  b_{n} ( \bar{x} ) \; \; \forall n \in
     \mathbb{N} , \forall  \bar{x} \in \Sigma^{\infty}
\end{equation}
Hence:
\begin{equation}
    r \in [ a_{n} ( \bar{x} ) , b_{n} ( \bar{x} ) ] \; \; \forall
    r \in \mathbb{R} , \forall n \in \mathbb{N} , \forall \bar{x}
    \in \Sigma^{\infty}
\end{equation}
so that:
\begin{equation}
       | S_{Dedekind} ( a_{0} , b_{0} ; \bar{x} ) | \; = \; \max (
       | \mathbb{R} | , |
    \mathbb{H}_{infinitesimals} |) \; = \;  |
    \mathbb{H}_{infinitesimals} |
\end{equation}
where we have used the definition \ref{def:hypercontinuity} and
the fact that the axiom \ref{ax:Generalized Continuum Hypothesis}
implies that:
\begin{equation}
    | \mathbb{R} | \; = \; \aleph_{1}
\end{equation}
\end{proof}

\smallskip

Considering together all the different cases it follows that:

\begin{proposition} \label{prop:hypercontinuous hyperreals are non-dedekindians 1}
\end{proposition}

\begin{hypothesis}
\end{hypothesis}
\begin{center}
 $ \mathbb{H} $ is hypercontinuous
\end{center}

\begin{thesis}
\end{thesis}

\begin{equation}
  \mathbb{H} \; = \; \mathbb{G}_{|
   \mathbb{H}_{infinitesimals} |}
\end{equation}
\begin{proof}
 The thesis is an immediate consequence of Proposition \ref{prop:first
 case}, Proposition \ref{prop:second case} and Proposition \ref{prop:third
 case}.
\end{proof}

\smallskip

Proposition \ref{prop:hypercontinuous hyperreals are
non-dedekindians 1} implies that:
\begin{proposition} \label{prop:hypercontinuous hyperreal numbers are non-Dedekindian}
\end{proposition}
\begin{hypothesis}
\end{hypothesis}
\begin{center}
  $ \mathbb{H} $ is hypercontinuous
\end{center}

\begin{thesis}
\end{thesis}
\begin{equation}
  \mathbb{H} \; \subset \; \mathbb{ND}
\end{equation}

Let us now consider the particular hyperreal number system
 $   ^{\star} \mathbb{R} $ of Nonstandard Analysis introduced in the definition
\ref{def:hyperreal number system of Non-Standard Analysis}.

Let us observe first of all that:
\begin{proposition} \label{prop:the particular hyperreal number system of Non-Standard analysis is hypercontinuous}
\end{proposition}
\begin{center}
 $ ^{\star} \mathbb{R} $ is hypercontinuous
\end{center}
\begin{proof}
  Demanding to \cite{Schechter-98}, \cite{Robinson-96},
\cite{Goldblatt-98} for all the mathematical-logical  details let
us recall that:
\begin{enumerate}
    \item there exists a map, called the \emph{$ \star$-tranform},
associating to each sentence $ \phi $ of the language $
\mathcal{L}_{\mathcal{R}} $ formalizing $ \mathbb{R} $ a sentence
$^{ \star} \phi $ of the language $ \mathcal{L}_{ ^{\star}
\mathcal{R}} $ formalizing $ ^{\star} \mathbb{R} $

\item there exists a principle, called the \emph{Transfer
Principle}, stating that a $ \mathcal{L}_{\mathcal{R}}$-sentence $
\phi $ is true if and only if $^{ \star } \phi $ is true
\end{enumerate}

Let us then consider the following $
\mathcal{L}_{\mathcal{R}}$-sentence:
\begin{equation}
    \phi \; := \; \forall x_{1},x_{2} \in \mathbb{R} \; \exists f:
    [ x_{1} , x_{2} ] \mapsto \mathbb{R} \; \text{ bijective}
\end{equation}
Applying to $ \phi $ the $ \star$-transform we obtain the
following $ \mathcal{L}_{ ^{\star} \mathcal{R}} $-sentence:
\begin{equation}
    ^{\star} \phi \; := \; \forall x_{1},x_{2} \in \: ^{\star} \mathbb{R} \; \exists f:
    [ x_{1} , x_{2} ] \mapsto  \, ^{\star}  \mathbb{R} \; \text{ bijective}
\end{equation}
Since $ \phi $ is true, the application of the Transfer Principle
allows to infer that $  ^{\star} \phi $ is true.

Choosing in particular $ x_{1} := + \, \epsilon \in \: ^{\star}
\mathbb{R}_{infinitesimals} $ and $ x_{2} := - \,  \epsilon \in \,
^{\star} \mathbb{R}_{infinitesimals} $ it follows that:
\begin{equation}
    | [ - \, \epsilon , \epsilon ] | \; = \; | ^{\star}
    \mathbb{R} | \; \geq \; \aleph_{1} \; \; \forall  \epsilon \in \,  ^{\star}
\mathbb{R}_{infinitesimals}
\end{equation}
from which it follows that:
\begin{equation}
    | ^{\star} \mathbb{R}_{infinitesimals} | \; = \; |   | [ - \, \epsilon , \epsilon ]
    | \; \geq \; \aleph_{1} \; \; \forall  \epsilon \in \,  ^{\star}
\mathbb{R}_{infinitesimals}
\end{equation}
Since, given $ \epsilon_{1}, \epsilon_{2} , \epsilon_{3} ,
\epsilon_{4} \in \, ^{\star} \mathbb{R}_{infinitesimals} \, : \,
\epsilon_{1} < \epsilon_{2} < \epsilon_{3} < \epsilon_{4}  $,
there exists clearly a bijection $ g : [ \epsilon_{1},
\epsilon_{2} ] \mapsto [ \epsilon_{3} , \epsilon_{4} ] $, it
follows that all the infinitesimal intervals have the same
cardinality; hence:
\begin{equation}
    | [ \epsilon_{1} , \epsilon_{2} ] | \; = \; |  ^{\star} \mathbb{R}_{infinitesimals}
    | \; \geq \; \aleph_{1} \; \; \forall \epsilon_{1},
    \epsilon_{2} \in \,  ^{\star} \mathbb{R}_{infinitesimals} \; :
    \; \epsilon_{1} < \epsilon_{2}
\end{equation}
\end{proof}

We can then infer that:
\begin{proposition} \label{prop:nonstandard analysis versus non-dedekindian numbers}
\end{proposition}
\begin{equation}
    ^{\star} \mathbb{R} \; \subset \; \mathbb{ND}
\end{equation}
\begin{proof}
  The thesis follows applying the Proposition \ref{prop:the particular hyperreal number system of Non-Standard analysis is
  hypercontinuous} and the Proposition \ref{prop:hypercontinuous hyperreal numbers are non-Dedekindian}
\end{proof}
\newpage
\section{Algorithmically random non-Dedekindian numbers}
Given $ n \in \mathbb{N}_{+} $ and $ a_{0}, b_{0} \in
\mathbb{G}_{n} $ such that $ a_{0} < b_{0} $:
\begin{definition} \label{def:algorithmically random generalized numbers}
\end{definition}
\emph{set of the random generalized numbers of order n with
respect to $  [ a_{0} , b_{0} ] $:}
\begin{equation}
    RANDOM ( \mathbb{G}_{n} ; [ a_{0} , b_{0} ]) \; := \;  \{ d \in S_{Dedekind} ( a_{0} , b_{0} ; \bar{x} ) :  \bar{x} \in RANDOM ( \Sigma^{\infty}) \}
\end{equation}

The definition\ref{def:algorithmically random generalized numbers}
is a generalization to non-Dedekindian numbers of the notion of
Martin L\"{o}f-Solovay-Chaitin algorithmic randomness as it is
shown by the following:

\begin{proposition}
\end{proposition}
\begin{equation}
    RANDOM ( \mathbb{G}_{1} ; [ a_{0} , b_{0} ]) \; = \; RANDOM (
    [ a_{0} , b_{0} ] )
\end{equation}
\begin{proof}
The thesis follows applying the Proposition \ref{prop:generalized
number of first order are the reals}, the definition
\ref{def:algorithmically random generalized numbers} and the
definition \ref{def:algorithmically random numbers of  an
interval}.
\end{proof}

\smallskip

As a particular case of the definition \ref{def:algorithmically
random generalized numbers} we have extended the notion of
algorithmic-randomness to the particular hyperreal number system
$^{\star} \mathbb{R} $ of Nonstandard Analysis introduced in the
definition \ref{def:hyperreal number system of Non-Standard
Analysis}.

Actually:

\begin{proposition}
\end{proposition}
\begin{equation}
    RANDOM ( ^{\star} \mathbb{R} ; [ 0 , 1] ) \; = \; \{ x \, + \, \epsilon \; : \; x \in RANDOM(0,1) , \epsilon \in \;  ^{\star} \mathbb{R}_{infinitesimals}   \}
\end{equation}
\begin{proof}
The thesis follows by the definition \ref{def:algorithmically
random generalized numbers} and the Proposition
\ref{prop:nonstandard analysis versus non-dedekindian numbers}
\end{proof}
\newpage
\appendix
\section{Chain ordered fields} \label{sec:Chain ordered fields}
Let us recall the following:
\begin{definition}
\end{definition}
\emph{field:}

a triple $ ( F , + , \cdot ) $ where:
\begin{itemize}
    \item F is a non-empty set
    \item $ + : F \times F \mapsto F , \cdot : F \times F \mapsto
    F $ are maps satisfying the following conditions:
    \begin{enumerate}
        \item commutativity of the sum:
\begin{equation}
    a + b \; = \; b + a \; \; \forall a,b \in F
\end{equation}
        \item associativity of the sum
\begin{equation}
    ( a + b ) + c \; = \; a + ( b+c) \; \;  \forall a,b,c \in F
\end{equation}
        \item existence  of the  zero element with respect to the the sum
\begin{equation}
    \exists 0 \in F \; : \; ( a + 0 \, = \, a \; \; \forall a \in
    F )
\end{equation}
        \item existence of the opposites with respect to the sum
\begin{equation}
    \forall a \in F \, \exists \, - a \in F \; : \; a + ( - a ) =
    0
\end{equation}
        \item commutativity of the product
\begin{equation}
    a \cdot b \; = \; b \cdot a \; \; \forall a,b \in F
\end{equation}
        \item associativity of the product
\begin{equation}
     ( a \cdot b ) \cdot c \; = \; a \cdot ( b \cdot c) \; \;  \forall a,b,c \in F
\end{equation}
        \item existence of the unity with respect to the product
\begin{equation}
    \exists 1 \in F \; : 1 \neq 0 \; \wedge \; ( 1 \cdot a \, =
    a \cdot 1 \, = \; a   \; \; \forall a \in F )
\end{equation}
        \item existence of the inverse with respect to the product
\begin{equation}
    \forall  a \in F : a \neq 0 \; \exists a^{-1} \in F \; : \; a \cdot
    a^{-1} \, = \, a^{-1} \cdot a \, = \, 1
\end{equation}
    \item distributivity of the product with respect to the sum
\begin{equation}
    a \cdot ( b + c) \; = \; a \cdot b + a \cdot c \; \; \forall a
    , b ,c \in F
\end{equation}
    \end{enumerate}
\end{itemize}

\smallskip

Let us recall that given a set S:
\begin{definition}
\end{definition}
\emph{partial ordering on S:}

$ \preceq \in \mathcal{P} ( S \times S )$ such that:
\begin{enumerate}
    \item reflectivity
\begin{equation}
    a \preceq a \; \; \forall a \in S
\end{equation}
    \item transitivity
\begin{equation}
  [( a \preceq b \; \wedge \; b \preceq c ) \; \Rightarrow \; a \preceq
  c ] \; \; \forall a,b,c \in S
\end{equation}
    \item antisimmetricity
\begin{equation}
  [ ( a \preceq b \, \wedge \, b \preceq a) \; \Rightarrow \; a =
  b ] \; \; \forall a ,b  \in S
\end{equation}
\end{enumerate}

\begin{definition}
\end{definition}
\emph{total ordering on S:}

$ \preceq $ partial ordering on S such that:
\begin{equation}
   ( a \preceq b \; \vee b \preceq a ) \; \; \forall a ,b \in S
\end{equation}

Given a field F:

\begin{definition}
\end{definition}
\emph{chain ordering on F:}

$ \preceq $ total ordering on F such that:
\begin{enumerate}
    \item translation invariance of the ordering
\begin{equation}
    (a \preceq b \; \Rightarrow \; a + c \preceq b + c ) \; \;
    \forall a , b , c \in F
\end{equation}
    \item positive elements are closed under product
\begin{equation}
  [( 0 \preceq a \; \wedge \; 0 \preceq b ) \; \Rightarrow \; 0 \preceq a
  \cdot b ]  \; \;
    \forall a , b  \in F
\end{equation}
\end{enumerate}

\begin{definition}
\end{definition}
\emph{chain ordered field:}

a couple $ ( F , \preceq ) $ such that:
\begin{enumerate}
    \item F is a field
    \item $ \preceq $ is a chain-ordering over F
\end{enumerate}

\smallskip

Given a chain ordered field  $ ( F , \preceq ) $ and $a, b \in F$:
\begin{definition}
\end{definition}
\begin{equation}
    a \prec b \; := \; a \preceq b \, \wedge \, a \neq b
\end{equation}

Given $ a , b \in F : a \prec b $
\begin{definition}
\end{definition}
\begin{equation}
    [ a , b ] \; := \; \{ x \in F \, : \, a \preceq x \preceq b \}
\end{equation}

Given $ S \subseteq F : S \neq \emptyset $:
\begin{definition}
\end{definition}
\emph{upper bounds of S:}
\begin{equation}
    UB(S) \; := \{ b \in F : a \preceq b \; \; \forall a \in S \}
\end{equation}
\begin{definition} \label{def:Dedekind completeness}
\end{definition}
\emph{$ ( F , \preceq ) $ is Dedekind complete:}
\begin{equation}
    \forall S \subseteq F : S \neq \emptyset \; \wedge \;  UB(S) \neq
    \emptyset \; \exists \sup(S) \in UB(S) \, : \, (  \sup(S)
    \preceq b \; \; \forall b \in UB(S) )
\end{equation}

\begin{definition}
\end{definition}
\emph{natural action of $ \mathbb{Z}$ over F:}

$ \circ : \mathbb{Z} \times F \mapsto F$:
\begin{equation}
    n \circ a \; := \; \sum_{i=1}^{n} a \; \;  n \in \mathbb{N} , a \in F
\end{equation}
\begin{equation}
    ( - n ) \circ a \; := \; - ( n \circ a ) \; \; n \in \mathbb{N} , a \in F
\end{equation}

\begin{definition} \label{def:Archimedean}
\end{definition}
\emph{ $ ( F , \preceq ) $ is Archimedean:}
\begin{equation}
   ( \exists b \in F \, : \;  n \circ a \preceq b \; \forall n \in
   \mathbb{Z} )
    \; \Rightarrow \; a = 0
\end{equation}
\newpage
\section{Martin Lof-Solovay-Chaitin randomness} \label{sec:algorithmic randomness}
In this section we will briefly review the definition of
algorithmically-random binary sequences.

Given a number $ n \in {\mathbb{N}} : n \geq 2 $  let us
introduce, preliminarily, the following:
\begin{definition}
\end{definition}
\emph{n-ary alphabet:}
\begin{equation}
    \Sigma_{n} \; := \; \{ k \in \mathbb{N} \, : \, k \leq n-1 \}
\end{equation}
Obviously:
\begin{proposition}
\end{proposition}
\begin{equation}
     \Sigma_{2} \; = \; \Sigma
\end{equation}

Denoted by $ \Sigma_{n}^{\star} := \cup_{k \in \mathbb{N}_{+}}
\Sigma_{n}^{k} $ the set of all the n-ary strings and by $
\Sigma_{n}^{\infty} := \{ \bar{x} : \mathbb{N}_{+} \mapsto
\Sigma_{n} \} $ the set of all the n-ary sequences, let us
introduce the following:
\begin{definition} \label{def:n-ary value}
\end{definition}
\emph{n-ary value}:

the map $ v_{n} : \Sigma_{n}^{\infty} \mapsto [ 0 , 1]$:
\begin{equation}
 v_{n} ( \bar{x} ) \; := \; \sum_{i=1}^{\infty } \frac{x_{i}}{n^{i}}
\end{equation}
and the more usual notation:
\begin{equation}
  (0.x_{1} \cdots x_{m} \cdots
)_{n} \; := \; v_{n} ( \bar{x} ) \; \; \bar{x} \in
\Sigma_{n}^{\infty}
\end{equation}

Let us introduce furthermore the following:
\begin{definition}
\end{definition}
\emph{n-ary nonterminating natural positional representation}:

the map $ r_{n} : [ 0 , 1] \mapsto \Sigma_{n}^{\infty} $:
\begin{equation}
  r_{n} ( (0.x_{1} \cdots x_{i} \cdots
)_{n} ) \;  := \; \bar{x}
\end{equation}
with the nonterminating condition requiring that the numbers of
the form $ (0. x_{1} \cdots x_{i} \overline{(n-1)})_{n} \, = \,
(0. \cdots (x_{i}+1) \bar{0})_{n}$ are mapped into the sequence $
x_{1} \cdots x_{i} \overline{(n-1)} $.

\smallskip

 Given $ n_{1} , n_{2} \in {\mathbb{N}} : \min ( n_{1} ,
n_{2} ) \geq 2 $:
\begin{definition}
\end{definition}
\emph{change of basis from $ n_{1} $ to $ n_{2} $}:

the map $ cb_{n_{1},n_{2}} : \Sigma_{n_{1}}^{\infty} \, \mapsto \,
\Sigma_{n_{2}}^{\infty} $:
\begin{equation}
  cb_{n_{1},n_{2}} ( \bar{x} ) \; := \; r_{n_{2}}(v_{n_{1}} (
  \bar{x}))
\end{equation}

Given $ X \in \Sigma_{n}^{\star} $:
\begin{definition}
\end{definition}
\begin{equation}
    X  \Sigma_{n}^{\infty} \; := \; \{ \bar{x} \in \Sigma_{n}^{\infty}  \, : ( \exists n \in \mathbb{N}_{+} \,: \, \vec{x} (n) \in S)   \}
\end{equation}

\smallskip

Endowed $ \Sigma_{n} $ with the discrete topology and $
\Sigma_{n}^{\infty} $ with the induced product topology $ \tau $:
\begin{definition}
\end{definition}
\emph{$ G \subset \Sigma_{n}^{\infty}$ is a constructively-open
set:}
\begin{enumerate}
    \item
\begin{equation}
    G \in \tau
\end{equation}
    \item
\begin{equation}
    \exists X \subset \Sigma_{n}^{\star} \text{ recursively enumerable
    } \, : \, G = X \Sigma_{n}^{\infty}
\end{equation}
\end{enumerate}
where we demand to \cite{Odifreddi-89} for the definition of
recursive enumerability.

\smallskip

\begin{definition}
\end{definition}
\emph{constructive sequence of constructively open sets (c.s.c.o.
sets):}

$ \{ G_{k} \, , \, k \in \mathbb{N} : k \geq 1 \} $ sequence of
constructively-open sets $ G_{k} \, = \, X_{k} \Sigma_{n}^{\infty}
$ such that:
\begin{equation}
    \exists \, X \subset \Sigma_{n}^{\star} \times \mathbb{N} \text{ recursively enumerable
    } \, : \;  X_{k} = \{ \vec{x} \in \Sigma_{n}^{\star} \, : \,
    ( \vec{x} , k) \in X \}
\end{equation}

\smallskip

\begin{definition}
\end{definition}
\emph{cylinder set with respect to $ \vec{x}  \, \in \,
\Sigma_{n}^{\star} $}:
\begin{equation} \label{eq:cylinder set}
\Gamma_{\vec{x}} \; := \; \{ \bar{y}  \in \Sigma_{n}^{\infty} \; :
\; \vec{y}(|\vec{x}|) = \vec{x} \}
\end{equation}

\begin{definition}
\end{definition}
\emph{cylinder - $ \sigma $ - algebra on $ \Sigma_{n}^{\infty} $}:
\begin{equation}
  {\mathcal{F}}_{cylinder} \; := \;  \sigma- \text{algebra generated by}  \{ \Gamma_{\vec{x}} \, : \, \vec{x} \in \Sigma_{n}^{\star}
  \}
\end{equation}

\begin{definition}
\end{definition}
\emph{Lebesgue measure:}

the probability measure over the measurable space $ (
\Sigma_{n}^{\infty} ,\mathcal{F}_{cylinder}  ) $:
\begin{equation}
      \mu_{Lebesgue} ( \Gamma_{\vec{x}} ) \; := \;  \frac{1}{n^{| \vec{x} |}} \; \;
  \forall \, \vec{x} \, \in \, \Sigma^{\star}
\end{equation}

\smallskip

Given $ S \subset \Sigma_{n}^{\infty} $:
\begin{definition}
\end{definition}
\emph{S is a constructively null set:}
\begin{multline}
    \exists \, \{ G_{k} \, , \, k \in \mathbb{N} : k \geq 1 \}
    \text{ c.s.c.o. sets  } \; : \; S \subset \cap_{k \geq 1}
    G_{k} \; \wedge \; \lim_{k \rightarrow + \infty}
    \mu_{Lebesgue} ( G_{k} ) \, = \, 0 \; constructively
\end{multline}

We can finally introduce the following:
\begin{definition}
\end{definition}
\emph{Martin L\"{o}f -  Solovay - Chatin random sequences over $
\Sigma_{n}$:}
\begin{equation}
    RANDOM( \Sigma_{n}^{\infty} ) \; := \; \Sigma_{n}^{\infty} - \{ S
    \subset \Sigma_{n}^{\infty}
    \, \text{ constructively null set } \}
\end{equation}

A key feature of the Martin L\"{o}f - Solovay - Chaitin notion of
algorithmic-randomness is the following \cite{Calude-02}:
\begin{proposition} \label{prop:basis-independence of randomness}
\end{proposition}
\emph{Basis-independence of randomness}:
\begin{equation}
  RANDOM( \Sigma_{n_{2}}^{\infty} ) \; = \; cb_{n_{1},n_{2}} (
  RANDOM(  \Sigma_{n_{1}}^{\infty} )) \; \; \forall n_{1} , n_{2} \in {\mathbb{N}} : \min ( n_{1} ,
n_{2} ) \geq 2
\end{equation}

 Proposition \ref{prop:basis-independence of randomness} allows
to restrict the analysis to algorithmically random binary
sequences without any lost of generality and to introduce the
following:

\begin{definition} \label{def:random numbers between zero and one}
\end{definition}
\emph{set of the algorithmically random numbers in the interval $
[ 0 , 1]$}:
\begin{equation}
    RANDOM(0,1) \; := \; \{ v_{2} ( \bar{x}) , \bar{x} \in
    RANDOM (\Sigma^{\infty}) \}
\end{equation}
\newpage
\section{The particular hyperreal number system of Non-Standard Analysis}
In this paper we work within the formal system ZFC, i.e. the
Zermelo-Fraenkel axiomatization of Set Theory augmented with the
Axiom of Choice (axiom \ref{def:Axiom of Choice}).

 Given a set $S \neq \emptyset $:
\begin{definition}
\end{definition}
\emph{filter on S:}
\begin{equation}
    \mathcal{F} \subseteq \mathcal{P} (S) \, : \, ( A \cap B \in
    \mathcal{F} \; \; \forall A , B \in \mathcal{F} ) \wedge ( A
    \in \mathcal{F} \wedge A \subseteq B \subseteq S \, \Rightarrow
    \, B \in \mathcal{F} )
\end{equation}
\begin{definition}
\end{definition}
\emph{ultrafilter on S:}

a filter $  \mathcal{F} $ on S such that:
\begin{equation}
   \mathcal{F} \neq \mathcal{P} (S) \; \wedge \; ( A \in
   \mathcal{F} \vee S - A \in \mathcal{F} \; \; \forall A \in  \mathcal{P}
   (S) )
\end{equation}

Given $ B \in \mathcal{P} (S) : B \neq \emptyset $:
\begin{definition}
\end{definition}
\emph{principal filter generated by  B:}
\begin{equation}
    \mathcal{F}^{B} \; := \; \{ A \in \mathcal{P} (S) \, : \, A \supseteq B \}
\end{equation}

In this paper we assume the following:

\begin{axiom} \label{def:Axiom of Choice}
\end{axiom}
\emph{Axiom of Choice:}
\begin{equation}
    \exists f : \mathcal{P} (S) \mapsto \cup_{B \in \mathcal{P}
    (S)} B \; : \; f(A) \in A \; \; \forall A \in  \mathcal{P} (S) :
    A \neq \emptyset
\end{equation}

A consequence of the axiom \ref{def:Axiom of Choice} is the
following \cite{Goldblatt-98}:

\begin{proposition} \label{prop:exists non principal ultrafilter over any infinite set}
\end{proposition}
\begin{equation}
    | S | \geq \aleph_{0} \; \Rightarrow \; \exists  \, \mathcal{F} \:
    \text {nonprincipal ultrafilter on S}
\end{equation}

Uniforming our notation to the one adopted for sequences over
finite alphabets let us introduce the following:
\begin{definition}
\end{definition}
\begin{equation}
    \mathbb{R}^{\infty} \; := \{ \bar{r} : \mathbb{N} \mapsto
    \mathbb{R} \}
\end{equation}

Given $ \bar{r} = \{ r_{n} \}_{n \in \mathbb{N}} ,  \bar{s} = \{
s_{n} \}_{n \in \mathbb{N}} \in \mathbb{R}^{\infty}$:

\begin{definition}
\end{definition}
\begin{equation}
    \bar{r} \oplus \bar{s} \; := \;  \{ r_{n} + s_{n} \}_{n \in \mathbb{N}}
\end{equation}
\begin{equation}
    \bar{r} \odot \bar{s} \; := \;  \{ r_{n} \cdot s_{n} \}_{n \in \mathbb{N}}
\end{equation}

\newpage
Given $ x \in \mathbb{R} $:
\begin{definition}
\end{definition}
\begin{equation}
    x^{\infty} \; := \; \{ r_{n} \}_{n \in \mathbb{N}} \in
    \mathbb{R}^{\infty} \, : \, r_{n} = x \; \; \forall n \in
    \mathbb{N}
\end{equation}

Let us now introduce the following:
\begin{definition}
\end{definition}
\begin{equation}
    NPU( \mathbb{N} ) \; := \; \{ \mathcal{F} \:
    \text {nonprincipal ultrafilter on } \mathbb{N} \}
\end{equation}

By Proposition \ref{prop:exists non principal ultrafilter over any
infinite set} it follows that:
\begin{proposition}
\end{proposition}
\begin{equation}
    NPU( \mathbb{N} )  \; \neq \; \emptyset
\end{equation}

Given $ \mathcal{F} \in  NPU( \mathbb{N} ) $ and $ \bar{r} = \{
r_{n} \}_{n \in \mathbb{N}} ,  \bar{s} = \{ s_{n} \}_{n \in
\mathbb{N}} \in \mathbb{R}^{\infty}$:
\begin{definition}
\end{definition}
\emph{$ \bar{r} $ and  $ \bar{s} $ are equal $ \mathcal{F}$-almost
everywhere:}
\begin{equation}
 \bar{r} \sim_{\mathcal{F}} \bar{s} \; := \; \{ n \in \mathbb{N} \, : \, r_{n} = s_{n}
 \} \in \mathcal{F}
\end{equation}

It may be proved that \cite{Goldblatt-98}:
\begin{proposition}
\end{proposition}
\begin{center}
  $ \sim_{\mathcal{F}} $ is an equivalence relation over $
  \mathbb{R}^{\infty}$
\end{center}

Let us finally introduce the following:
\begin{definition}
\end{definition}
\begin{equation}
    ^{\star} \mathbb{R}_{\mathcal{F}} \; = \; \frac{\mathbb{R}^{\infty}}{\sim_{\mathcal{F}}}
\end{equation}

Given $ \bar{r} = \{ r_{n} \}_{n \in \mathbb{N}} ,  \bar{s} = \{
s_{n} \}_{n \in \mathbb{N}} \in \mathbb{R}^{\infty}$:
\begin{definition}
\end{definition}
\begin{enumerate}
    \item
 \begin{equation}
    [ \bar{r} ]_{\mathcal{F}} + [ \bar{s} ]_{\mathcal{F}} \; := \;
    [ \bar{r} \oplus  \bar{s} ]_{\mathcal{F}}
\end{equation}
    \item
\begin{equation}
    [ \bar{r} ]_{\mathcal{F}} \cdot [ \bar{s} ]_{\mathcal{F}} \; := \;
    [ \bar{r} \odot  \bar{s} ]_{\mathcal{F}}
\end{equation}
    \item
\begin{equation}
    [ \bar{r} ]_{\mathcal{F}} \leq [ \bar{s} ]_{\mathcal{F}} \; := \;
    \{ n \in \mathbb{N} \, : \, r_{n} \leq s_{n} \} \in {\mathcal{F}}
\end{equation}
\end{enumerate}
\newpage
It may be proved that \cite{Goldblatt-98}:
\begin{proposition}
\end{proposition}
\begin{center}
  $ (  ^{\star} \mathbb{R}_{\mathcal{F}} , + , \cdot , \leq ) $ is
  an hyperreal number system $ \forall \mathcal{F} \in  NPU( \mathbb{N}
  ) $ with zero $ [ 0^{\infty} ]_{\mathcal{F}} $ and unity $ [ 1^{\infty}
  ]_{\mathcal{F}} $ and where $  x \in \mathbb{R} $ is identified with $ [ x^{\infty} ]_{\mathcal{F}} $
\end{center}
where we have used the definition \ref{def:hyperreal number
system} of an hyperreal number system.

Furthermore the assumption of the axiom \ref{ax:Generalized
Continuum Hypothesis} implies that:
\begin{proposition} \label{prop:independence from the nonprincipal ultrafilter}
\end{proposition}
\begin{center}
   $ (  ^{\star} \mathbb{R}_{\mathcal{F}_{1}} , + , \cdot , \leq )
   $ is isomorphic to $ (  ^{\star} \mathbb{R}_{\mathcal{F}_{2}} , + , \cdot , \leq
   )  \; \; \forall \mathcal{F}_{1} ,  \mathcal{F}_{2} \in NPU( \mathbb{N}
  ) $
\end{center}

Proposition \ref{prop:independence from the nonprincipal
ultrafilter} allows to give the following:

\begin{definition} \label{def:hyperreal number system of Non-Standard Analysis}
\end{definition}
\emph{hyperreal number system of Non-Standard Analysis:}
\begin{equation}
    (  ^{\star} \mathbb{R} , + , \cdot , \leq ) \; := \; ( ^{\star} \mathbb{R}_{\mathcal{F}} , + , \cdot , \leq
    ) \; \;  \mathcal{F} \in NPU( \mathbb{N})
\end{equation}

\newpage
\section{Notation}
\begin{center}
  \begin{tabular}{|c|c|}
  \hline
  $ i.e. $ & id est \\
  $ \forall $ & for all (universal quantificator) \\
  $ \exists $ & exists (existential quantificator) \\
  $ \exists \, !  $ & exists and is unique \\
   $ x \; = \; y $ & x is equal to y \\
  $ x \; := \; y $ & x is defined as y \\
  $ \wedge $ & and (logical conjunction) \\
   $ \vee $ & or (logical conjunction) \\
  $ gcd( n, m) $  & greatest common divisor of n and m \\
  $ \Sigma $ & binary alphabet \\
  $ \Sigma^{\star} $ & set of the binary strings \\
  $ \Sigma^{\infty} $ & set of the binary sequences \\
  $ RANDOM( \Sigma^{\infty}) $ & set of the random  binary sequences  \\
  $ \vec{x} $ & binary string \\
  $ \bar{x} $ & binary sequence \\
   $ x_{n} $ & $ n^{th} $ digit of the string $\vec{x} $ or of the sequence $ \bar{x} $ \\
    $ \vec{x}(n) $ & prefix of length n of the string $ \vec{x} $ or of the sequence $ \bar{x} $  \\
    $ | \vec{x} | $ & length of the string $ \vec{x} $ \\
    $ \vec{x} ^{\infty}  $ & sequence made of infinite repetitions
of the string $ \vec{x} $ \\
   $ \cdot $ & concatenation operator \\
  $ \lambda $ & empty string \\
  $| S |$ & cardinality of the set $ S $ \\
  $ \mathcal{P} (S)$ & power set of the set S \\
  $ \aleph_{n} $ & $ n^{th} $ infinite cardinal number \\
  $ \mathbb{N} $ & set of the natural numbers \\
  $ \omega $  & ordinality of  $ \mathbb{N} $ \\
   $ \mathbb{N}_{+} $ & set of the strictly positive natural numbers \\
 $ \mathbb{Z} $ & set of the integer numbers \\
   $ \mathbb{Q} $ & set of the rational numbers \\
  $ \mathbb{R} $ & set of the real numbers \\
  $ \mathbb{H} $ & an hyperreal number system \\
  $ ^{\star} \mathbb{R}$ & the hyperreal number system of
  Non-standard Analysis \\
  $  \mathbb{H}_{infinitesimals} $ & set of the infinitesimal
  elements of $  \mathbb{H} $ \\
  $  \mathbb{H}_{unboundeds} $ & set of the unbounded
  elements of $  \mathbb{H} $ \\
  $ std(x)$ & standard part of x \\
  $ \mathcal{D}_{a_{0},b_{0}} $ & Dedekind operator with respect to $
  [ a_{0} , b_{0} ] $ \\
  $ RANDOM ( [a_{0} , b_{0}] ) $ & real random numbers with respect to $ [a_{0} , b_{0}] $ \\
  $ S_{Dedekind} ( a_{0} , b_{0} ; \bar{x} ) $ & Dedekind set with respect to $ [ a_{0} , b_{0}] $ and
$\bar{x}$ \\
  $ \mathbb{G}_{n} $ & set of the generalized numbers of order n \\
  $ \mathbb{N} \mathbb{D} $ & set of the non-Dedekindian numbers \\
 $ RANDOM(  \mathbb{G}_{n} ; [ a_{0} , b_{0} ] )$ & set of the random generalized numbers of order
 n with respect to $[a_{0}, b_{0} ]$ \\
  \hline
\end{tabular}
\end{center}
\newpage
\section{Acknoledgements}
I would like to thank prof. Cristian Calude for many precious
remarks and suggestions.

\end{document}